\documentclass[11pt]{article}

\usepackage[margin=1.3in]{geometry}
\usepackage[scaled]{helvet}  

\usepackage[english]{babel}
\usepackage[utf8]{inputenc}
\usepackage{amsmath}
\usepackage{amsfonts,dsfont}
\usepackage{amsthm}
\usepackage{graphicx}
\usepackage{program,color}
\usepackage{natbib}
\usepackage{algorithmicx}
\usepackage[ruled]{algorithm}
\usepackage{algpseudocode}
\usepackage{caption}
\usepackage{subcaption}

\newtheorem{definition}{Definition}
\newtheorem{theorem}{Theorem}

\newtheorem{lemma}{Lemma}
\newtheorem{proposition}{Proposition}

\newcommand{\MMD}{\mathrm{MMD}}
\newcommand{\ED}{\mathrm{ED}}

\newcommand{\E}{\mathbb{E}}

\newcommand{\Hh}{\mathcal{H}}
\newcommand{\RR}{\mathbb{R}}

\newcommand{\esssup}{ess sup}
\newcommand{\supp}{supp}

\newcommand{\one}{{\mathbf{1}}} 

\newcommand{\dotprod}[2]{\ensuremath{\langle #1 , #2\,\rangle}} 
\DeclareMathOperator{\tr}{trace}
\title{ On Wasserstein Two Sample Testing and Related Families of Nonparametric Tests
} 

\author{
Aaditya Ramdas$^*$\\
University of California, Berkeley\\
\texttt{aramdas@cs.berkeley.edu}
\and
Nicol\'as Garc\'ia Trillos$^*$\\
Brown University\\
\texttt{ngarciat@andrew.cmu.edu}
\and
Marco Cuturi\\
Kyoto University\\
\texttt{mcuturi@i.kyoto-u.ac.jp}
}

\date{\today}

\begin{document}
\maketitle

\begin{abstract}
Nonparametric two sample or homogeneity testing is a decision theoretic problem that involves identifying differences between two random variables without making parametric assumptions about their underlying distributions. The literature is old and rich, with a wide variety of statistics having being intelligently designed and analyzed, both for the unidimensional and the multivariate setting. Our contribution is to tie together many of these tests, drawing connections between seemingly very different statistics. In this work, our central object is the Wasserstein distance, as we form a chain of connections from univariate methods like the Kolmogorov-Smirnov test, PP/QQ plots and ROC/ODC curves, to multivariate tests involving energy statistics and kernel based maximum mean discrepancy. Some connections proceed through the construction of a \textit{smoothed} Wasserstein distance, and others through the pursuit of a ``distribution-free'' Wasserstein test. Some observations in this chain are implicit in the literature, while others seem to have not been noticed thus far. Given nonparametric two sample testing's classical and continued importance, we aim to provide useful connections for theorists and practitioners familiar with one subset of methods but not others. 
\end{abstract}

\section{Introduction}

Nonparametric two sample testing (or homogeneity testing) deals with detecting differences between two $d$-dimensional distributions, given samples from both, without making any parametric distributional assumptions. The popular tests for $d=1$ are rather different from those for $d> 1$, and our interest is in tying together different tests used in both settings.
There is a massive literature on the two-sample problem, having been formally studied for nearly a century, and there is no way we can cover the breadth of this huge and historic body of work. Our aim is much more restricted --- we wish to study this problem through the eyes of the beautiful Wasserstein distance. We wish to form connections between several seemingly distinct families of such tests, both intuitively and formally, in the hope of informing both practitioners and theorists who may have familiarity with some sets of tests, but not others. We will also only introduce related work that has a direct relationship with this paper.

There are also a large number of tests for \textit{parametric} two-sample testing (assuming a form for underlying distributions, like Gaussianity), and yet others for testing only differences in \textit{means} of distributions (like Hotelling's t-test, Wilcoxon's signed rank test, Mood's median test). Our focus will be much more restricted --- in this paper, we will restrict our attention only to \textit{nonparametric} tests for testing differences in (any moment of the underlying) \textit{distribution}.

Our paper started as an attempt to understand testing with the Wasserstein distance (also called earth-mover's distance or transportation distance). The main prior work in this area involve studying the ``trimmed'' comparison of distributions by \cite{munk1998nonparametric,freitag2007nonparametric} with applications to biostatistics, specifically population bioequivalence, and later by \cite{esteban08,esteban12}. Apart from two-sample testing, the study of univariate \textit{goodness-of-fit testing} (or one-sample testing) was undertaken in \cite{delbarrio99,delbarrio00,delbarrio05}, and summarized extremely well in \cite{delbarrio04notes}. There are other semiparametric  works specific to goodness-of-fit testing for location-scale families that we do not mention here, since they diverge from our interest in fully nonparametric two-sample testing for generic distributions.  

In this paper, we uncover an interesting relationship  between the multivariate Wasserstein test and the (Euclidean) Energy distance test, also called the Cramer test, proposed independently by \cite{energydistance} and \cite{baringhausfranz04}. This proceeds through the construction of a \textit{smoothed Wasserstein distance}, by adding an entropic penalty/regularization --- varying the weight of the regularization interpolates between the Wasserstein distance at one extreme and the Energy distance at the other extreme. This also gives rise to a new connection between the univariate Wasserstein test and popular univariate data analysis tools like quantile-quantile (QQ) plots and the Cramer von-Mises (CvM) test. Due to the relationship between distances and kernels, we will also establish connections to the kernel-based multivariate test by \cite{mmd} called the Maximum Mean Discrepancy, or MMD. Finally, the desire to design a univariate \textit{distribution-free} Wasserstein test will lead us to the formal study of Receiver Operating Characteristic (ROC) curves, relating to work by \cite{roc96}.

Intuitively, the underlying reasons for the similarities and differences between these above tests can be seen through two lenses. First is the \textit{population} viewpoint of how different tests work with different \textit{representations} of distributions; most of these tests are based on differences between  quantities that completely specify a distribution --- (a) cumulative distribution functions (CDFs), (b) quantile functions (QFs), and (c) characteristic functions (CFs). Second is from the \textit{sample} viewpoint of the behavior these statistics show under the null hypothesis; most of these tests have null distributions based on norms of Brownian bridges, alternatively viewed as infinite sums of weighted chi-squared distributions (due to the Karhunen-Loeve expansion).

While we connect a wide variety of popular and seemingly disparate families of tests, there are still further classes of tests that we do not have space to discuss. Some examples of tests quite different from the ones studied here include rank based tests as covered by the excellent book \cite{lehmann06ranks}, and graphical tests that include spanning tree methods by \cite{friedmanrafsky79} (generalizing the runs test by \cite{waldwolfowitz40}), nearest-neighbor based tests by \cite{schilling86} and \cite{henze88}, and the cross-match tests by \cite{rosenbaum}. The book by \cite{thas2010} is also a very useful reference.

\paragraph{Paper Outline and Contributions.} The rest of this paper proceeds as follows. In Section \ref{sec:NTST}, we formally present the notation and setup of nonparametric two sample testing, as well as briefly introduce three different ways of comparing distributions ---using CDFs, QFs and CFs.  In Section \ref{sec:WTSTmulti} we  form a novel connection between the multivariate Wasserstein distance, to the multivariate Energy Distance, and to the kernel MMD, through an entropy-smoothed Wasserstein distance. In Section \ref{sec:WTSTuni} we relate the univariate Wasserstein two-sample test to PP and QQ plots/tests. Lastly, in Section \ref{sec:WTSTdf}, we will design a  univariate Wasserstein test statistic that is also ``distribution-free'' unlike its classical counterpart, providing a careful and rigorous analysis of its limiting distribution by connecting it to ROC/ODC curves.

\section{Nonparametric Two Sample Testing}\label{sec:NTST}

 More formally, given i.i.d. samples $X_1,...,X_n \sim P$ and $Y_1,...,Y_m \sim Q$, where $P$ and $Q$ are probability measures on $\mathbb{R}^d$.  We denote by $P_n$ and $Q_m$ the corresponding empirical measures. A test  $\eta $ is a function from the data $D_{m,n} := \{X_1,...X_n,Y_1,...,Y_m\} \in \RR^{d(m+n)}$ to $\{0,1\}$ (or to $[0,1]$ if it is a randomized test). 

Most tests proceed by calculating a scalar test statistic $T_{m,n} := T(D_{m,n}) \in \RR$ and deciding $H_0$ or $H_1$ depending on whether $T_{m,n}$, after suitable normalization, is smaller or larger than a threshold $t_\alpha$. $t_{\alpha}$ is calculated based on a prespecified false positive rate $\alpha$, chosen so that, $\E_{H_0} \eta \leq \alpha$, at least asymptotically. Indeed, all tests considered in this paper are of the form
$$
\eta(X_1,...,X_n, Y_1,...,Y_m) = \mathbb{I}\left( T_{m,n} > t_\alpha \right)
$$
We follow the Neyman-Pearson paradigm, were a test is judged by its power  $\phi = \phi(m,n,d,P,Q,\alpha) = \E_{H_1} \eta$. We say that a test $\eta$ is consistent, in the classical sense, when
$$
\phi \rightarrow 1 \mbox{ as } m,n \rightarrow \infty, \alpha \to 0.
$$

All the tests we consider in this paper will be consistent in the classical sense mentioned above. Establishing general conditions under which these tests are consistent in the high-dimensional setting is largely open.
All the test statistics considered here are of the form that they are typically small under $H_0$ and large under $H_1$ (usually with appropriate scaling, they converge to zero and to infinity respectively with infinite samples). The aforementioned threshold $t_\alpha$ will be determined by the distribution of the test statistic being used under the null hypothesis (i.e. assuming the null was true, we would like to know the typical variation of the statistic, and we reject the null if our observation is far from what is typically expected under the null). This naturally leads us to study the \textit{null distribution} of our test statistic, i.e. the distribution of our statistic under the null hypothesis. Since these are crucial to running and understanding the corresponding tests, we will pursue their description in some detail in this paper.

\subsection{Three Ways to Compare Distributions}

The literature broadly has three dominant ways of comparing distributions, both in one and in multiple dimensions. These are based on three different ways of characterizing distributions --- cumulative distribution functions (CDFs), characteristic functions (CFs) and quantile functions (QFs). Many of the tests we will consider involve calculating differences between (empirical estimates of) these quantities.

For example, it is well known that the Kolmogorov-Smirnov (KS) test by \cite{kolmogorov33} and \cite{smirnov48} involves differences in empirical CDFs. We shall later see that in one dimension, the Wasserstein distance calculates differences in QFs.

The KS test, the related Cramer von-Mises criterion by \cite{cramer28} and \cite{vonmises28}, and Anderson-Darling test by \cite{anderson1952asymptotic} are very popular in one dimension, but their usage has been  more restricted in higher dimensions. This is mostly due to the curse of dimensionality involved with estimating multivariate empirical CDFs. While there has been  work on generalizing these popular one-dimensional to higher dimensions, like \cite{bickel69}, these are seemingly not the most common multivariate tests.

Two classes of tests that are actually quite popular are kernel and distance based tests. As we will recap in more detail in later sections, it is  known that the Gaussian kernel MMD implicitly calculates a (weighted) difference in CFs and the Euclidean energy distance implicitly works with a difference in (projected) CDFs.


\section{Entropy Smoothed Wasserstein Distances}\label{sec:WTSTmulti}
The theory of optimal transport (see \citep{villani09}) provides a set of powerful tools to compare probability measures and distributions on $\mathbb{R}^d$ through the knowledge of a metric on $\mathbb{R}^{d}$, which we assume to be the usual Euclidean metric in what follows. Among that set of tools, the following family of $p$-Wasserstein distances between probability measures is the best known.

\subsection{Wasserstein Distance}\label{subsec:wass}
Given an exponent $p\geq 1$, the definition of the $p$-Wasserstein distance reads:

\begin{definition}[Wasserstein Distances]\label{def:was}
For $p\in[1,\infty)$ and Borel probability measures $P,Q$ on $\RR^d$ with finite $p$-moments, their $p$-Wasserstein distance \citep[Sect. 6]{villani09} is
\begin{equation}\label{eq:wass}
	W_p(P,Q)=\left(\inf_{\pi\in\Gamma(P,Q)}\int_{\RR^d \times \RR^d} \| X-Y \|^p d\pi \right)^{1/p},
\end{equation}
where $\Gamma(P,Q)$ is the set of all joint probability measures on $\mathbb{R}^d\times \mathbb{R}^d$ whose marginals are $P,Q$, i.e. such that for all subsets $A \subset \mathbb{R}^d$ we have $\pi(A\times \mathbb{R}^d)=P(A)$ and $\pi(\mathbb{R}^d\times A)=Q(A)$.
\end{definition}

A remarkable feature of Wasserstein distances is that Definition~\ref{def:was} applies to all measures regardless of their absolute continuity with respect to the Lebesgue measure: the same definition works for both empirical measures and for their densities if they exist.

Writing $\one_n$ for the $n$-dimensional vector of ones, when comparing two empirical measures
with uniform\footnote{The Wasserstein machinery works also for non-uniform weights. We do not mention this in this paper because all of the measures we consider in the context of two-sample testing are uniform.} weight vectors $\mathbf{1}_n/n$ and $\mathbf{1}_m/m$,
the Wasserstein distance $W_p(P_n,Q_m)$ exponentiated to the power $p$ is the optimum of a network flow problem known as the transportation problem \citep[Section 7.2]{bertsimas1997introduction}. This problem has a linear objective and a polyhedral feasible set, defined respectively through the matrix $M_{XY}$ of pairwise distances between elements of $X$ and $Y$ raised to the power $p$,
\begin{equation}\label{eq:distMatrix}M_{XY}:= [\| X_i-Y_j \|^p]_{ij} \in\mathbb{R}^{n\times m},\end{equation} 
and the polytope $U_{nm}$ defined as the set of $n\times m$ nonnegative matrices such that their row and column marginals are equal to $\mathbf{1}_n/n$ and $\mathbf{1}_m/m$ respectively:
\begin{equation}\label{eq:polytope} 
	U_{nm}:= \{ T\in\mathbb{R}_+^{n\times m}\;:\; T \mathbf{1}_m = \mathbf{1}_n/n,\, T^T \mathbf{1}_n = \mathbf{1}_m/m \}.
\end{equation}

Let $\dotprod{A}{B}:= \tr(A^T B)$ be the usual Frobenius dot-product of matrices.  Combining Eq.~\eqref{eq:distMatrix} and \eqref{eq:polytope}, we have that $W_p^p(P_n,Q_m)$ is the optimum of a linear program $S$ of $n\times m$ variables,
\begin{equation}\label{eq:primal}
W_p^p(P_n,Q_m) =\min_{T\in U_{nm}}\dotprod{T}{M_{XY}},
\end{equation}
of feasible set $U_{nm}$ and cost matrix $M_{XY}$.

We finish this section by pointing out that  the rate of convergence as $n,m \rightarrow \infty$ of $W_p(P_n,Q_m)$ towards $W_p(P,Q)$ gets slower as the dimension $d$ grows under mild assumptions. For simplicity of exposition consider $m=n$. For any $p \in [1,\infty)$, it follows from \cite{Dud} that for $d\geq 3$, the difference between $W_p(P_n, Q_n)$ and $W_p(P,Q)$ scales as $n^{-1/d}$. We also point out that when $d=2$ the rate actually scales as $\frac{\sqrt{\ln(n)}}{\sqrt{n}}$ (see \cite{AKT}).  Finally, we note that when considering $p=\infty$ the rates of convergence are different to those when $1\leq p <\infty$ . The work of \cite{LeightonShor,ShorYukich,garcia14}  show that the rate of convergence of $W_\infty(P_n,Q_n)$ towards $W_\infty(P,Q)$ is of order $\left(\frac{\ln(n)}{n} \right)^{1/d}$ when $d\geq 3$ and $\frac{(\ln(n))^{3/4}}{n^{1/2}}$ when $d=2$. Hence, the original Wasserstein distance by itself may not be a favorable choice for a multivariate two sample test.

\subsection{Smoothed Wasserstein Distance}
Aside from the slow convergence rate of the Wasserstein distance between samples from two different measures to their distance in population, computing the optimum of \eqref{eq:primal} is expensive. This can be easily seen by noticing that the transportation problem boils down to an optimal assignment problem when $n=m$. Since the resolution of the latter has a cubic cost in $n$, all known algorithms that can solve the optimal transport problem scale at least super-cubicly in $n$. Using an idea that can be traced back as far as \citet{Schrodinger31}, \citet{cuturi2013sinkhorn}  recently proposed to use an entropic regularization of the optimal transport problem, to define the Sinkhorn divergence between  $P,Q$ parameterized by $\lambda\geq 0$ as

\begin{equation}\label{eq:sink}S_\lambda^p(P,Q) := \min_{T\in U_{nm}}\,\lambda\dotprod{T}{M_{XY}}-E(T).
\end{equation}
where $E(T)$ is the entropy of $T$ seen as a discrete joint probability distribution, namely $E(T):=-\sum_{ij}T_{ij} \log(T_{ij}).$
Let $T_\lambda$ be the minimizer of the above smoothed optimal transport problem.
 
This approach has two benefits: \emph{(i)} because $E(T)$ is $1$-strongly convex with respect to the $\ell_1$ norm, the regularized problem is itself strongly convex and admits a unique optimal solution $T_\lambda$, as opposed to the initial OT problem, for which the minimizer may not be unique; \emph{(ii)} this optimal solution $T_\lambda$ is a diagonal scaling of $e^{-M_{XY}}$, the element-wise exponential matrix of $-M_{XY}$. One can easily show using the Lagrange method of multipliers that there must exist two non-negative vectors $u\in\mathbb{R}^n, v\in\mathbb{R}^m$ such that $T_\lambda := D_u e^{-M_{XY}} D_v$, where $D_u$ $D_v$ are diagonal matrices with $u$ and $v$ on their diagonal. The solution to this diagonal scaling problem can be found efficiently through Sinkhorn's algorithm \citep{Sinkhorn67}, which has a linear convergence rate \citep{franklin1989scaling}. Sinkhorn's algorithm can be implemented in a few lines of code that only require matrix vector products and elementary operations,  hence easily parallelized on modern hardware.


\subsection{Smoothing the Wasserstein Distance to Energy Distance}

An interesting class of  tests are distance-based ``energy statistics'' as introduced in parallel by \cite{baringhausfranz04} and  \cite{energydistance}.  
The statistic, called the \textit{Cramer statistic} by the former paper and \textit{Energy Distance} by the latter,  corresponds to the population quantity
$$
\ED ~:=~  2\E \|X-Y\| - \E \|X-X'\|  - \E\|Y-Y'\|, 
$$
where  $X,X' \sim P$ and $Y,Y' \sim Q$ (all i.i.d.).
An associated test statistic can be calculated as 
\begin{eqnarray}
\ED_b &:=&   
\frac{2}{mn}\sum_{i=1}^n\sum_{j=1}^m \|X_i - Y_j\| - \frac1{n^2}\sum_{i, j=1}^n \|X_i - X_j\| - \frac1{m^2}\sum_{i, j=1}^m \|Y_i - Y_j\| \nonumber .
\label{eq:hd}
\end{eqnarray}

Remarkably, $\ED(P,Q) = 0$ iff $P=Q$. Hence, rejecting when $\ED_b$ is larger than an appropriate threshold leads to a test which is consistent against all fixed alternatives where $P \neq Q$ under mild conditions (like finiteness of $\E[X],\E[Y]$); see aforementioned references for details. 
Then, the Sinkhorn divergence defined in ~\eqref{eq:sink} can be linked to the the energy distance when the regularization parameter is set to $\lambda=0$, through the following formula:
\begin{equation}
\ED_b = 2S_0^1(P_n,Q_m) - S_0^1(P_n,P_n) - S_0^1(Q_m,Q_m).
\label{eqnEDSinkhorn}
\end{equation}
Indeed, notice first that $T_0$ is the maximal entropy table in $U_{nm}$, namely the outer product $(\one_n \one_m^T)/nm$ of the marginals $\one_n/n$ and $\one_m/m$. Then \eqref{eqnEDSinkhorn} follows from the observation
$$S_0^1(P_n,Q_m) = \frac{1}{nm} \sum_{i,j} \| X_i-Y_j \| , \quad S_0^1(P_n,P_n) = \frac{1}{n^2} \sum_{i,j=1}^n \| X_i-X_j \|, \quad S_0^1(Q_m,Q_m) = \frac{1}{m^2} \sum_{i,j=1}^m \| Y_i-Y_j \|.$$ 

\subsection{From Energy Distance to Kernel Maximum Mean Discrepancy}

Another popular class of tests that has emerged over the last decade, are kernel-based tests introduced independently by  \cite{mmd06} and \cite{alba08}, and expanded on in \cite{mmd}. Without getting into technicalities that are irrelevant for this paper, the \textit{Maximum Mean Discrepancy} between $P,Q$ is defined as
$$
\MMD(H_k,P,Q) := \max_{\|f\|_{\Hh_k} \leq 1} \E_P f(X) - \E_Q f(Y)
$$
where $\Hh_k$ is a Reproducing Kernel Hilbert Space associated with Mercer kernel $k(\cdot,\cdot)$, and $\|f\|_{\Hh_k} \leq 1$ is its unit norm ball. While it is easy to see that $\MMD \geq 0$ always, and also that $P=Q$ implies $\MMD=0$, \cite{mmd06} show that if $k$ is ``characteristic'', the equality holds iff $P=Q$ (the Gaussian kernel $k(a,b) = \exp(-\|a-b\|^2/\gamma^2)$ is a popular example). 
Using the Riesz representation theorem and the reproducing property of $\Hh_k$, one can argue that
$
\MMD(\Hh_k, P,Q) = \|\E_P k(X,.) - \E_Q k(Y,.)\|_{\Hh_k}
$
and hence using the reproducing property again, one can conclude that
$$
\MMD^2 = \E k(X,X') + \E k(Y,Y') - 2\E k(X,Y).
$$

This gives rise to a natural associated test statistic, a plugin estimator of $\MMD^2$:
\begin{eqnarray*}
\MMD^2_u(k(\cdot , \cdot)) &:=& \frac1{n^2}\sum_{i, j=1}^n k(X_i,X_j) +  \frac1{m^2}\sum_{i, j=1}^m k(Y_i,Y_j) - \frac{2}{mn}\sum_{i=1}^n \sum_{j=1}^m k(X_i,Y_j). \label{eq:hk}
\end{eqnarray*}
Apart from the fact that $\MMD(P,Q) = 0$ iff $P=Q$ the other fact that makes this a useful test statistic is that its estimation error, i.e. the error of $\MMD^2_u$ in estimating $\MMD^2$, scales like $\sqrt{\frac{m+n}{mn}}$, independent of $d$ \footnote{This is unlike KL-divergence, which is also zero iff $P=Q$, but is in general hard to estimate in high dimensions.}. See \cite{mmd} for a detailed proof. 

At first sight, the Energy Distance and the MMD look like fairly different tests. However, there is a natural connection that proceeds in two steps. Firstly, there is no reason to stick to only the Euclidean norm $\|\cdot \|_2$ to measure distances for ED --- the test can be extended to other norms, and in fact also other metrics; \cite{lyons} explains the details for the closely related independence testing problem. Following that,  \cite{hsiceqdcov} discuss the relationship between distances and kernels (again for independence testing, but the same arguments hold in the two sample testing setting also). Loosely speaking, for every kernel $k$, there exists a metric $d$ (and also vice versa), given by $d(x,y) := (k(x,x) + k(y,y))/2 - k(x,y)$, such that MMD with kernel $k$ equals ED with  metric $d$. This is a very strong connection between these two families of tests.

\section{Wasserstein Distance and PP or QQ tests}\label{sec:WTSTuni}

 For univariate random variables, a PP plot is a graphical way to view differences in empirical CDFs, while QQ plots are analogously for comparing QFs. Instead of relying on graphs, we can also make such tests more formal and rigorous as follows. We first present some results on the asymptotic distribution of the difference between $P_n$ and $Q_m$ when using the distance between the CDFs $F_n$ and $G_m$ and then later when using the distance between the QFs $F_n^{-1}$ and $G_m^{-1}$.
For simplicity we assume that both distributions $P$ and $Q$ are supported on the interval $[0,1]$; we remark that under mild assumptions on $P$ and $Q$, the results we present in this section still hold without such a boundedness assumption. Moreover we assume for simplicity that the CDFs $F$ and $G$ have positive densities on $[0,1]$.

\subsection{Comparing CDFs (PP)}

We start by noting $F_n$ may be interpreted as a random element taking values in the space $\mathcal{D}([0,1])$ of right continuous functions with left limits. It is well known that
\begin{equation}
\sqrt{n}\left(  F_n   - F \right) \rightarrow_{w}   \mathbb{B}\circ F 
\label{weakconvCDF}
\end{equation}
where $\mathbb{B}$ is a standard Brownian bridge in $[0,1]$ and where the weak convergence  $\rightarrow_{w}$ is understood as convergence of probability measures in the space $\mathcal{D}([0,1])$; see Chapter 3 in \cite{Billingsley} for details.  From this fact and the independence of the samples, it follows that under the null hypothesis $H_0:P=Q$, as $n, m \rightarrow \infty$
\begin{equation}
\sqrt{\frac{mn}{m+n}}  \left(   F_n - G_m    \right)  = \sqrt{\frac{mn}{m+n}} \left(F_n - F \right)  + \sqrt {\frac{mn}{m+n}} \left( G- G_m  \right)  \rightarrow_{w}   \mathbb{B}\circ F .
\label{weakConv}
\end{equation}

The previous fact, and continuity of the function $h \in \mathcal{D}([0,1]) \mapsto  \int_{0}^{1}(h(t))^2 dt$, imply that as $n, m \rightarrow \infty$, we have under the null,
$$   \frac{mn}{m+n}\int_{0}^{1}   \left(   F_n(t) - G_m(t)   \right)^2 dt  \rightarrow_{w} \int_{0}^{1} (\mathbb{B}(F(t)))^2 dt. $$
Observe that the above asymptotic null distribution 
depends on $F$  which is unknown in practice. This is an obstacle\footnote{This obstacle can be avoided in the goodness-of-fit testing context, when we want to test if $X_1,\dots, X_n$ was drawn from a \textit{known} CDF $F$ or not. This is the original purpose of  the $L^2$-statistics of the von Mises type. Briefly, from \eqref{weakconvCDF}, and since the function $f \in \mathcal{D}([0,1])  \mapsto \int_{0}^{1} (f(t))^2 d F(t)$ is continuous, we deduce
$$ \int_{0}^{1} ( F_n(t) - F(t) )^2 d F(t) \rightarrow_{w} \int_{0}^{1} \mathbb{B}(F(t))^2 dF(t) = \int_{0}^{1} (\mathbb{B}(s))^2 ds, $$
where the second equality follows by a change of variables, leading to an expression that does not depend on $F$.} when considering any $L^p $-distance, with $1 \leq p < \infty$, between the empirical cdfs $F_n$ and $G_m$. 
Luckily, a different situation occurs when one considers the $L^\infty$-distance between $F_n$ and $G_m$.  Under the null, using again  \eqref{weakconvCDF} we deduce that
\begin{equation}
  \sqrt{\frac{mn}{m+n}} \| F_n - G_m\|_\infty \rightarrow_w     \| \mathbb{B}\circ F\|_\infty=  \| \mathbb{B}\|_\infty,
\end{equation}
where the equality in the previous expression follows from the fact that the continuity of $F$ implies that the interval $[0,1]$ is mapped onto the interval $[0,1]$. 
This test statistic, the so-called Kolmogorov-Smirnov test statistic, is hence appropriate for two sample problems.

\subsection{Comparing QFs (QQ)}

We now turn our attention to $QQ$ (quantile-quantile) plots and specifically the $L^2$-distance between $F_n^{-1}$ and $G_m^{-1}$. It can be shown that if $F$ has a differentiable density $f$ which (for the sake of simplicity) we assume is bounded away from zero, then 
$$   \sqrt{n} (F_n^{-1}-F^{-1}) \rightarrow_w  \frac{\mathbb{B}}{f \circ F^{-1}}. $$
For a proof of the above statement see Chapter 18 in \cite{ShorackWellner}; for an alternative proof where the weak convergence is considered in the space of probability measures on $L^2((0,1))$ (as opposed to the space $\mathcal{D}([0,1])$ we have been considering thus far) see \cite{delbarrio04notes}.

We note that from the previous result and independence, it follows that under the null hypothesis $H_0: P=Q$,
$$   \sqrt{\frac{mn}{m+n}} (F_n^{-1}-G_m^{-1}) \rightarrow_w   \frac{\mathbb{B}}{f \circ F^{-1}} .$$
In particular by continuity of the function $h \in L^2((0,1)) \mapsto \int_0^{1} (h(t))^2 dt$, we deduce that
$$   \frac{mn}{m+n} \int_{0}^{1} (F_n^{-1}-G_m^{-1})^2 dt \rightarrow_w   \int_{0}^{1}  \frac{(\mathbb{B} (t))^2}{(f \circ F^{-1}(t))^2} dt .$$
Hence, as was the case when we considered  the difference of the cdfs $F_n$ and $G_m$, the asymptotic distribution of the $L^2$-difference (or analogously any $L^p$-difference for finite $p$) of the empirical quantile functions is also distribution dependent. Note however that there is an important difference between QQ and PP plots when using the $L^\infty$ norm. We saw that the asymptotic distribution of the $L^\infty$ norm of the difference of $F_n $ and $G_m$ is (under the null hypothesis) distribution free. Unfortunately, in the quantile case, we obtain
$$\sqrt{\frac{mn}{m+n}}\| F_n^{-1}- G_n^{-1}\|_{\infty}  \rightarrow_w  \| \frac{\mathbb{B}}{f \circ F^{-1}} \|_\infty,$$
which of course is distribution dependent. Since one would have to resort to computer-intensive Monte-Carlo techniques (like bootstrap or permutation testing) to control type-1 error, these tests are sometimes overlooked (though with modern computing speeds, they merit further study).

\subsection{Wasserstein is a QQ test}

We recall that in general, for $p\in[1,\infty)$ the $p$-Wasserstein distance between two probability measures $P, Q$ on $\RR$ with finite $p$-moments is given by
\begin{equation}
 W_p(P, Q ) := \inf_{\pi\in \Gamma(P,Q)}  \left( \int_{\RR \times \RR} \| x- y \|^p d \pi(x,y)   \right)^{1/p}. 
 \label{pWasser}
\end{equation}

Because the Wasserstein distance measures the cost of transporting mass from the original distribution $P$ into the target distribution $Q$, one can say that it measures "horizontal" discrepancies between $P$ and $Q$. Intuitively, two probability distributions $P$ and $Q$ that are different over "long" (horizontal) regions will be far away from each other in the Wasserstein distance sense, because in that case mass has to travel long distances to go from the original distribution to the target distribution.  In the one dimensional case (in contrast with what happens in dimension $d>1$), the $p$-Wasserstein distance has a simple interpretation in terms of the quantile functions $F^{-1}$ and $G^{-1}$ of $P$ and $Q$ respectively. The reason for this is that the optimal way to transport mass from $P$ to $Q$  has to satisfy certain monotonicity property which we describe in the proof of the following Lemma. This is a well known fact that can be found, for example, in \cite{thas2010}. Nevertheless here we present its proof in the Appendix for the sake of completeness. 

\begin{proposition}\label{prop:p-wass}
The $p$-Wasserstein distance between two probability measures $P$ and $Q$ on $\RR$ with $p$-finite moments can be written as
$$W_p^p(P, Q)   = \int_{0}^{1} \lvert F^{-1}(t)- G^{-1}(t)   \rvert^p dt  , $$
where $F^{-1}$ and $G^{-1}$ are the quantile functions of $P$ and $Q$ respectively.
\end{proposition}

Having considered the $p$-Wasserstein distance $W_p(P,Q)$ for $p \in [1,\infty)$ in the one dimensional case, we conclude this section by considering the case $p=\infty$.  Let $P, Q$ be two probability measures on $\RR$ with bounded support. That is, assume that there exist a number $N>0$ such that $\supp(P) \subseteq [-N,N]$ and $\supp(Q) \subseteq [-N,N]$. We define the $\infty$-Wasserstein distance between $P$ and $Q$ by
$$  W_\infty(P,Q):= \inf_{\pi \in \Gamma(P, Q)}  \esssup_{\pi} \lvert x-y \rvert. $$ 
Proceeding as in the case $p\in [1,\infty)$, it is possible to show that the $\infty$-Wasserstein distance between $P$ and $Q$ with bounded supports can be written in terms of the difference of the corresponding quantile functions as
$$ W_\infty(P,Q) = \| F^{-1} - G^{-1} \|_{\infty}.$$
 
 The Wasserstein distance is also sometimes called the Kantorovich-Rubinstein metric and the Mallow's distance in the statistical literature, where it has been studied extensively due to its ability to capture weak convergence precisely -- $W_p(F_n,F)$ converges to 0 if and only if $F_n$ converges in distribution to $F$ and also the $p$-th moment of $X$ under $F_n$ converges to the corresponding moment under $F$; see \cite{dobrushin1970prescribing,mallows1972note,bickel1981some}. 

\section{A Distribution-Free Wasserstein Test}\label{sec:WTSTdf}
As we earlier saw, under the null hypothesis $H_0: P=Q$, the statistic 
$ \frac{mn}{m+n} \int_{0}^{1} \left( F_n^{-1}(t)  - G_m^{-1}(t)\right)^2 d t  $
has an asymptotic distribution which is not distribution free, i.e., it depends on $F$. We also saw that as opposed to what happens with the asymptotic distribution of the $L^\infty$ distance between $F_n$ and $G_m$, the asymptotic distribution of $\|  F_n^{-1} - G_m^{-1} \|_\infty$ \textit{does} depend on the cdf $F$. 

In this section we show how we can construct a distribution-free Wasserstein test. To prove that it is distribution-free, we connect it to the theory of ROC and ODC curves. 

\subsection{Relating Wasserstein Distance to ROC and ODC curves}

Let $P$ and $Q$ be two distributions on $\RR$ with cdfs $F$ and $G$ and quantile functions $F^{-1}$ and $G^{-1}$ respectively. We define the \textit{ROC} curve between $F $ and $G$ as the function.
$$ ROC(t) := 1- F (G^{-1}(1-t)) ,\quad t \in [0,1].$$
In addition, we define their \textit{ODC} curve by,
$$ ODC(t) := G (F^{-1}(t)), \quad t \in [0,1]. $$
The following are straightforward  properties of the ROC curve (see \cite{Turnbull}). 
\begin{enumerate}
\item The $ROC$ curve is increasing and $ROC(0)=0$, $ROC(1)=1$.
\item If $G(t) \geq F(t)$ for all $t$ then $ROC(t) \geq t$ for all $t$.
\item If $F$ and $G$ have densities with monotone likelihood ratio, then the ROC curve is concave.
\item The area under the ROC curve is equal to  $\mathbb{P}(Y < X )$, where $Y \sim Q$ and $X \sim P$.
\end{enumerate}
Intuitively speaking, the faster the ROC curve increases towards the value $1$, the easier it is to distinguish the distributions $P$ and $Q$. Observe from their definitions that the ROC curve can be obtained from the ODC curve after reversing the axes.  Given this, we focus from this point on only one of them, the ODC curve being more convenient.

The first observation about the ODC curve is that it can be regarded as the quantile function of the distribution $G _{\sharp } P$ (the push forward of $P$ by $G$) on $[0,1]$ which is defined by
$$  G _\sharp P([0,\alpha)) :=  P \left(  G^{-1}\left( [ 0,\alpha) \right) \right), \quad \alpha \in [0,1]. $$
Similarly,  we can consider the measure ${G_{m}}_{\sharp} P_n  $, that is, the push forward of $P_n$ by $G_m$. We crucially note that the empirical ODC curve $G_m \circ F_n^{-1}$ is the quantile function of ${G_{m}}_{\sharp} P_n  $. 
From Section \ref{sec:WTSTuni}, we deduce that
$$  W_p^p( {G_{m}}_{\sharp} P_n, {G}_{\sharp} P  )  = \int_{0}^{1} \lvert   G_m\circ F_n^{-1}(t) - G \circ F^{-1}(t)  \rvert^p dt  $$
for every $p \in[1,\infty)$ and also
$$   W_\infty(  {G_{m}}_{\sharp} P_n, {G}_{\sharp} P  ) = \lvert \lvert   G_m\circ F_n^{-1} - G \circ F^{-1}     \rvert\rvert_\infty. $$
That is, the $p$-Wasserstein distance between the measures ${G_{m}}_{\sharp} P_n$ and $ {G}_{\sharp} P$ can be computed by considering the $L^p$ distance of the ODC curve and its empirical version.


%

%


First we argue that under the null hypothesis $H_0:P=Q$, the distribution of empirical ODC curve is actually independent of $P$. In particular,  $W_p^p( {G_{m}}_{\sharp} P_n, {G}_{\sharp} P  )$ and $W_\infty( {G_{m}}_{\sharp} P_n, {G}_{\sharp} P  )$ are distribution free under the null! This is the content of the next lemma, proved in the Appendix.
\begin{lemma}[Reduction to uniform distribution] Let $F,G$ be two continuous and strictly increasing CDFs and let $X_1,\dots, X_n \sim F$ and $Y_1,\dots, Y_m \sim G$ be two independent samples. We let $F_n$ and $G_m$ be the CDFs associated to the empirical distributions induced by the $X$s and the $Y$s respectively. Consider the (unknown) random variables, which are distributed uniformly on $[0,1]$,
$$   U^X_k := F( X_k ), \    U^Y_k := G(Y_k). $$
Let $F^U_n$ be the empirical CDF associated to the (uniform) $U^X$s and let $G^U_m$ be the empirical CDF associated to the (uniform) $U^Y$s. Then, under the null hypothesis $H_0: F=G$ we have
$$    G_m(X_k)  =  G^U_m(U^X_k), \quad \forall k  \in \left\{1, \dots, n \right\}. $$
 In particular,
$$ G_m \circ F_n^{-1}(t) = G^U_m \circ {F^U_n}^{-1}(t), \quad \forall t\in [0,1].  $$
\label{ReducUniform}
\end{lemma}

Note that since $U^X_k, U^Y_k$ are obviously instantiations of uniformly distributed random variables, the RHS of the last equation only involves uniform r.v.s and hence, the distribution of $G_m \circ F_n^{-1}$ is independent of $F,G$ under the null. Now we are almost done -- this above lemma will imply that the Wasserstein distance between $G_m \circ F_n^{-1}$ and the uniform distribution  $U[0,1]$ (since $G \circ F^{-1}(t)=t=U^{-1}(t)=U(t)$ for $t \in [0,1]$ when $G=F$) also does not depend on $F,G$.

More formally, we establish a result related to the asymptotic distribution of $W_p^p( {G_{m}}_{\sharp} P_n, {G}_{\sharp} P  )$ and $W_\infty( {G_{m}}_{\sharp} P_n, {G}_{\sharp} P  )$. We do this by first considering the asymptotic distribution of the difference between the empirical ODC curve and the population ODC curve regarding both of them as elements in the space $\mathcal{D}([0,1])$.
This is the content of the following Theorem which follows directly from the work of \cite{Komlos} (see \cite{Turnbull}).

\begin{theorem}
\label{TheoremROC}
Suppose that $F$ and $G$ are two cdfs  with densities $f,g$ satisfying
$$   \frac{g(F^{-1}(t))}{f(F^{-1}(t))} \leq C,  $$
for all $t\in [0,1]$. Also, assume that 
$ \frac{n}{m} \rightarrow \lambda \in [0,\infty) $
as $n, m \rightarrow \infty$. Then, 
$$   \sqrt{\frac{mn}{m+n}}\left( G_m (F_n^{-1}(\cdot ))  - G(F^{-1}(\cdot ))  \right)   \rightarrow_w \sqrt{\frac{\lambda}{\lambda+1}} B_1(G \circ F^{-1}(\cdot))  +\sqrt{ \frac{1}{\lambda + 1}} \frac{g(F^{-1}(\cdot))}{f(F^{-1}(\cdot))} B_2(\cdot), $$
where $B_1$ and $B_2$ are two independent Brownian bridges and where the weak convergence must be interpreted as weak convergence in the space of probability measures on the space $\mathcal{D}([0,1])$.
\end{theorem}

As a corollary, under the null hypothesis $H_0: P = Q$ we obtain the following. Suppose that the CDF $F$ of $P$ is continuous and strictly increasing. Then, 
$$ \frac{mn}{m+n} W_2^2( {G_{m}}_{\sharp} P_n, {G}_{\sharp} P) =  \frac{mn}{m+n}\int_{0}^{1}( G_m (F_n^{-1}(t)) -t )^2 dt \rightarrow_w \int_{0}^{1}(\mathbb{B}(t))^2dt   $$
and
$$  \sqrt{\frac{mn}{m+n}} W_\infty( {G_{m}}_{\sharp} P_n, {G}_{\sharp} P) =  \sqrt{\frac{mn}{m+n}} \sup_{t \in [0,1]} \lvert G_m (F_n^{-1}(t)) -t \rvert \rightarrow_w   \sup_{t \in [0,1]}\lvert \mathbb{B}(t) \rvert.   $$
To see this, note that by Lemma \ref{ReducUniform} it suffices to consider $F(t) = t$ in $[0,1]$. In that case, the assumptions of Theorem \ref{TheoremROC} are satisfied and the result follows directly.\\

The takeaway message of this section is that instead of considering the Wasserstein distance between $F_m$ and $G_n$, whose null distribution depends on unknown $F$, one can instead consider the Wasserstein distance between $G_m(F_n^{-1})$ and the uniform distribution $U[0,1]$, since its null distribution is independent of $F$, i.e. we have constructed a distribution-free test.  

\section{Conclusion}

In this paper, we connect a wide variety of univariate and multivariate test statistics, with the central piece being the Wasserstein distance. The Wasserstein statistic is closely related to univariate tests like the Kolmogorov-Smirnov test, graphical QQ plots, and a distribution-free variant of the test is proposed by connecting it to ROC/ODC curves. Through entropic smoothing, the Wasserstein test is also related to the multivariate tests of Energy Distance and hence transitively to the Kernel Maximum Mean Discrepancy. 

We hope that this is a useful resource to connect the seemingly vastly different families of two sample tests, many of which can be analyzed under the two umbrellas of our paper -- whether they differentiate between CDFs, QFs or CFs, and what their null distributions look like. A comprehensive empirical survey is also of interest but out of our current scope.

\subsection*{Acknowledgments}

AR was supported by the grant NSF IIS-1247658.

\bibliography{mmd}
\bibliographystyle{plainnat}

\newpage
\appendix

\section{Proof of Proposition \ref{prop:p-wass}}

\begin{proof}
We first observe that the \textit{infimum} in the definition of $W_p(P, Q)$ can be replaced by \textit{minimum}, namely, there exists a transportation plan $\pi \in \Gamma(P, Q)$ that achieves the infimum in \eqref{pWasser}. This can be deduced in a straightforward way by noting that the expression $ \int_{\RR \times \RR} \lvert x- y \rvert^p d \pi(x,y)$ is linear in $\pi$ and that the set $\Gamma(P,Q)$ is compact in the sense of weak convergence of probability measures on $\RR \times \RR$. Let us denote by $\pi^*$ an element in $\Gamma(P,Q)$ realizing the minimum in \eqref{pWasser}. Let $(x_1, y_1) \in \supp(\pi^*) $ and $(x_2,y_2) \in \supp(\pi^*)$ (here $\supp(\pi^*)$ stands for the support of $\pi$) and suppose that $x_1< x_2$. We claim that the optimality of $\pi^*$ implies that $y_1 \leq y_2$. To see this, suppose for the sake of contradiction that this is not the case, that is, suppose that $y_2 < y_1$. We claim that in that case
\begin{equation}
\lvert x_1-y_2 \rvert^p  + \lvert x_2-y_1 \rvert^p < \lvert x_1-y_1 \rvert^p + \lvert x_2-y_2 \rvert^p . 
\label{auxWassQQ}
\end{equation}
Note that for $p=1$ this follows in a straightforward way. For the case $p>1$, first note that $x_1 < x_2 $ and $y_2 < y_1$ imply that there exists $t\in (0,1)$ such that $tx_1+ (1-t)y_1 = tx_2 +(1-t)y_2$. Now, note that
$$\lvert x_1 -  y_2\rvert = \lvert x_1 - \left(tx_1 + (1-t)y_1 \right) \rvert  + \lvert  \left(tx_1 + (1-t)y_1 \right) - y_2 \rvert $$
because the points $x_1$, $y_2$ and $tx_1 + (1-t)y_1$ all lie on the same line segment.
But then, using the fact that $tx_1 + (1-t)y_1 = tx_2 + (1-t)y_2$, we can rewrite the previous expression as
$$ \lvert x_1 - y_2\rvert  = (1-t)\lvert x_1- y_1\rvert  + t \lvert y_2 - x_2 \rvert.$$
Using the strict convexity of the function $t \mapsto t^p$ ( when $p>1$), we deduce that 
$$  \lvert x_1 - y_2\rvert^p < (1-t) \lvert x_1- y_1\rvert^p + t \lvert x_2- y_2\rvert^p.$$
In a similar fashion, we obtain
$$  \lvert x_2 - y_1\rvert^p < t \lvert x_1- y_1\rvert^p + (1-t) \lvert x_2- y_2\rvert^p.$$
Adding the previous two inequalities we obtain \eqref{auxWassQQ}. Note however that \eqref{auxWassQQ} contradicts the optimality of $\pi^*$, because it shows that $\pi^*$ is not \textit{cyclically monotone}, which essentially means that it is possible to rearrange the way mass is transported from $P$ to $Q$ by $\pi^*$ in order to reduce the transportation cost (it would be cheaper to send mass from $x_1$ to $y_2$ and from $x_2$ to $y_1$ than to send mass from $x_1$ to $y_1$ and from $x_2$ to $y_2$). Therefore, we conclude that if $(x_1, y_1) \in \supp(\pi^*)$, $(x_2, y_2) \in \pi^*$ and $x_1 < x_2$, then $y_1 \leq y_2$.  

Now, for $x \in \supp(P)$ and $y \in \supp(Q)$ we claim that $(x,y) \in \supp(\pi^*)$ if and only if $F(x)=G(y)$. To see this note that from the monotonicity property just established we deduce that $(x,y) \in \supp(\pi^*)$ if and only if $ \pi^*(\RR, (-\infty, y])  = \pi^*((-\infty, x]), (-\infty,y]  ) = \pi^*((-\infty, x], \RR)$. In turn, the fact that $\pi^* \in \Gamma(P,Q)$ implies that $\pi^*((-\infty, x], \RR) = F(x)$ and $\pi^*( \RR, (-\infty,y]) = G(y)$. From the previous relation we conclude that
$$  \int_{\RR \times \RR} \lvert x-y \rvert^p d \pi^*(x,y) = \int_{\supp(\pi^*)}\lvert x-y \rvert^p  d\pi^*(x,y) = \int_{0}^{1}\lvert F^{-1}(t) - G^{-1}(t) \rvert^p dt,  $$
as we wanted to show.
\end{proof}

\section{Proof of Lemma \ref{ReducUniform}}

\begin{proof}
We denote by $Y_{(1)} \leq \dots \leq Y_{(m)}$ the order statistic associated to the $Y$s.  For $k=1,\dots,m-1$ and $t \in (0,1)$, we have $G_m(t) = \frac{k}{m}$ if and only if $t \in [Y_{(k)}, Y_{(k+1)} )$ which holds if and only if $t \in [ F^{-1}(U^Y_{(k)}) , F^{-1}(U^Y_{(k+1)}) )$, which in turn is equivalent to $F(t) \in [ U^Y_{(k)} , U^Y_{(k+1)}  )$.   Thus, $G_m(t) = \frac{k}{m}$ if and only if $G^U_m (F(t)) = \frac{k}{m}$. From the previous observations we conclude that $G_m = G^U_m\circ F$.  Finally, since $X_k  = F^{-1}( U^X_k )$ we conclude that

$$   G_m(X_k) =  G^U_m\circ F \circ F^{-1}(U^X_k)= G^U_m( U^X_k).  $$
\end{proof}

\end{document}